\documentclass[amscd,amssymb,verbatim]{amsart} 
\input epsf.sty
\newcommand{\al}{{\alpha}}

\newcommand{\be}{\begin{enumerate}}
\newcommand{\bl}{{\rm{Bl}\,}}

\newcommand{\ca}{{\mathcal A}}
\newcommand{\cb}{{\mathcal B}}
\newcommand{\cc}{{\mathcal C}}

\newcommand{\cf}{{\mathcal F}}
\newcommand{\cg}{{\mathcal G}}
\newcommand{\ci}{{\mathcal I}}

\newcommand{\cl}{{\mathcal L}}
\newcommand{\cm}{{\mathcal M}}
\newcommand{\cn}{{\mathcal N}}
\newcommand{\co}{{\mathcal O}}

\newcommand{\cs}{{\mathcal S}}

\newcommand{\dc}{{\mathbb C}}
\newcommand{\dr}{{\mathbb R}}

\newcommand{\ee}{\end{enumerate}}

\newcommand{\ga}{\Gamma}

\newcommand{\lmt}{\longmapsto}
\newcommand{\lra}{\longrightarrow}

\newcommand{\pr}{\noindent{\bf Proof. }}

\newcommand{\ra}{\rightarrow}

\newcommand{\sm}{\setminus}
\newcommand{\Span}{\text{\rm span}}

\newcommand{\wti}{\widetilde }

\newtheorem{thm}{Theorem}[section]
\newtheorem{df}  [thm]{Definition}

\newtheorem{lm}  [thm]{Lemma}

\newtheorem{prop}[thm]{Proposition}

\newtheorem{rem} [thm]{Remark}
\newtheorem{expl} [thm]{Example}

\numberwithin{equation}{section}

\newenvironment{remrm}{\begin{rem} \rm}{\end{rem}}
\newenvironment{explrm}{\begin{expl} \rm}{\end{expl}}

\title 
       {Incidence combinatorics of resolutions}
 
\author{Eva-Maria Feichtner and Dmitry N. Kozlov}
                
\date{\noindent
November 2000, final version June 2003\\[0.05cm]
      \hskip15pt The second author acknowledges support by the 
Forschungsinstitut f\"ur Mathematik, ETH Z\"urich, 
and the Swedish National Science Foundation.\\[0.05cm]
 \hskip15pt
MSC 2000 Classification: 
  primary 06A07; 
secondary 52C35, 
          05E99, 
          14E15, 
          14M25. \\
Keywords: arrangement models, 
          resolution of singularities,
          blowups,
          building sets, 
          toric varieties. 
}
  
\address{Departement Mathematik, ETH Z\"urich, 8092 Z\"urich, Switzerland}
\email{feichtne@math.ethz.ch}

\address{Department of Mathematics, KTH Stockholm, 100 44 Stockholm, Sweden}
\email{kozlov@math.ethz.ch, kozlov@math.kth.se} 

\begin{document}

\begin{abstract}    
We introduce notions of combinatorial blowups,  building sets, and
nested sets for arbitrary meet-semilattices. This gives a~common abstract 
framework for the incidence combinatorics occurring in the context of 
De Concini-Procesi models of subspace arrangements and  resolutions 
of singularities in toric varieties.
Our main theorem states that a sequence of combinatorial blowups, 
prescribed by a building set in linear extension compatible order, 
gives the face poset of the corresponding simplicial complex of nested sets. 
As applications we trace the incidence combinatorics through every 
step of the De Concini-Procesi model construction, and we introduce
the notions of  building sets and nested sets to the context of toric 
varieties.

There are several other instances, such as models of stratified
manifolds and certain graded algebras associated with finite lattices, 
where our combinatorial framework has been put to work; we present 
an outline in the end of this paper.
\end{abstract}

\maketitle

              
\section{Introduction}\label{sect_intr}
   For an arbitrary meet-semilattice we introduce notions of
{\it combinatorial blowups, building sets\/}, and {\it nested sets}. 
The definitions are given on a~purely order-theoretic level without 
any reference to geometry. This provides a common abstract framework 
for the incidence combinatorics occurring in at least two different
situations in algebraic geometry: the construction of De Concini-Procesi 
models of subspace arrangements~\cite{DP1}, and the resolution of 
singularities in toric varieties.

  The various parts of this abstract framework have received 
different emphasis within different situations: while the notion of
combinatorial blowups clearly specializes to stellar subdivisions of 
defining fans in the context of toric varieties, building sets and 
nested sets were introduced in the context of model constructions 
by De Concini \& Procesi~\cite{DP1} (earlier and in a more special 
setting by Fulton~\& MacPherson~\cite{FM}), from where we adopt our 
terminology. This correspondence however is not complete: 
the building sets in \cite{DP1, FM} are not canonical, they depend 
on the geometry, while ours do not. See Section~\ref{ssect_DPmodels} 
for further details.

  It was proved in \cite{DP1} that a sequence of blowups within
an arrangement  of complex linear subspaces leads from
the~intersection stratification of complex space given by the 
maximal subspaces of the arrangement to an arrangement model 
stratified by divisors with normal crossings. 
In the context of toric varieties, there exist many different
procedures for stellar subdivisions of a defining fan that result
in a simplicial fan, so-called simplicial resolutions.  

  The purpose of our Main Theorem \ref{thm_main} is to unify 
these two situations on the combinatorial level: a sequence of 
combinatorial blowups, performed on a~(combinatorial) building 
set in linear extension compatible order, transforms the initial
semilattice to a semilattice  where all intervals are boolean 
algebras, more precisely to the face poset of the corresponding
simplicial complex of nested sets. In particular, the structure of 
the resulting semilattice can be fully described by the initial 
data of nested sets. Both the formulation and the proof of our 
main theorem are purely combinatorial.

We sketch the content of this article:

\noindent
{\bf Section 2.}
After providing some basic poset terminology, we define building sets
and nested sets for meet-semilattices in purely order-theoretic terms
and develop general structure theory for these notions.
 
\noindent  
{\bf Section 3.}
  We define combinatorial blowups of meet-semilattices, and study their
effect on building sets and nested sets. The section contains our 
Main Theorem~\ref{thm_main} which describes the result of blowing up 
the elements of a building set in terms of the initial nested set complex.
   
\noindent  
{\bf Section 4.} 
This section is devoted to relating our abstract framework to two different
contexts in algebraic geometry: 
In~\ref{ssect_DPmodels}
we briefly review the construction of De Concini-Procesi models for subspace 
arrangements. We show that the change of the incidence combinatorics of the 
stratification in a single construction step is described by a 
combinatorial blowup of the semilattice of strata. 
In~\ref{ssect_tv} we draw the connection to simplicial 
resolutions of toric varieties: we recognize stellar subdivisions as 
combinatorial blowups of the face posets of defining fans and discuss 
the notions of building and nested sets in this context.

\noindent  
{\bf Section 5.} Since a first version of this paper has been written
and has circulated in the fall of 2000, our combinatorial framework
for the incidence combinatorics of resolutions has been taken up in
various contexts. We outline the model construction for real subspace
and halfspace arrangements and for real stratified manifolds by
G.~Gaiffi~\cite{Ga2}. Moreover, we give a short account of the study
of a graded algebra associated with any finite lattice in~\cite{FY},
where our combinatorial generalization of originally geometric notions
leads to the construction of an, at first sight, unrelated geometric
counterpart for wonderful models of hyperplane arrangements. We also
note that, more recently, our combinatorial resolutions were studied
in the context of log resolutions of arrangement ideals~\cite{BS}.


\section{Building sets and nested sets of meet-semilattices} 
\label{sect_buidg_nest_sets}

\subsection{Poset Terminology} \label{ssect_pos_term}
\mbox{ } 

\noindent 
We recall some notions from the theory of partially ordered sets, and
refer to~\cite[Ch.\,3]{St} for further details. All posets discussed in 
this paper will be finite. A poset~$\cl$ is called a {\em meet-semilattice\/}
if any two elements $x,y\,{\in}\,\cl$ have a greatest lower bound, i.e.,
the set $\{z\,{\in}\,\cl\,|\,z\,{\leq}\,x, z\,{\leq}\,y\}$ has a maximal 
element,
called the {\em meet}, $x\,{\wedge}\,y$, of~$x$ and~$y$. Greatest lower
bounds of subsets $A\,{=}\,\{a_1,\ldots,a_t\}$ in~$\cl$ we denote with
$\bigwedge A\,{=}\, a_1\,{\wedge}\,\ldots \,{\wedge}\, a_t$. In particular, 
meet-semilattices have a unique minimal element denoted~$\hat 0$. Minimal
elements in $\cl\,{\setminus}\,\{\hat 0\}$ are called the {\em atoms} 
in~$\cl$. Meet-semilattices share the following property: for any subset
$A\,{=}\,\{a_1,\ldots,a_t\}\,{\subseteq}\,\cl$ the set 
$\{x\,{\in}\,\cl \,|\, x\,{\geq}\, a \text{ for all } a\,{\in}\, A\}$ 
is either empty 
or it has a unique minimal element, called the~{\em join}, $\bigvee A\,{=}\, 
a_1\,{\vee}\,\ldots \,{\vee}\, a_t$, of~$A$. If the meet-semilattice needs 
to be specified, we write $(\bigvee A)_{\cl}\,{=}\, 
(a_1\,{\vee}\,\ldots \,{\vee}\, a_t)_{\cl}$ for the join of~$A$ in~$\cl$.
For brevity, we talk about semilattices throughout the paper, meaning 
meet-semilattices. 

 Let $P$ be an arbitrary poset.
For $x\,{\in}\,P$ set:
$P_{\leq x}\,{=}\,\{y\,{\in}\,P\,{|}\, y\,{\leq}\,x\}$; $P_{<x}$, 
and $P_{\geq x}$, $P_{>x}$ are defined analogously. For subsets 
$\cg\,{\subseteq}\,P$ 
with the induced order, we define $\cg_{\leq x}\,{=}\,\{y\,{\in}\,P\,{|}\, 
y\,{\in}\,\cg, y\,{\leq}\,x\}$, and
$\cg_{<x}$  again analogously. For intervals in~$P$ we use the standard
notations $[x,y]\,{:=}\,\{z\,{\in}\, P\,|\, x\,{\leq}\,z \,{\leq}\,y\}$,
$[x,y)\,{:=}\,\{z\,{\in}\, P\,|\, x\,{\leq}\,z \,{<}\,y\}$, etc.

A poset is called {\it irreducible\/} if it is not a direct product of 
two other posets, both consisting of at least two elements. For a poset
$P$ with a unique minimal element~$\hat 0$,
we call $I(P)=\{x\,{\in}\,P\,|\,[\hat 0,x]\,\, \text{is 
irreducible}\,\}$ the {\em set of irreducible elements\/}
in~$P$.
In particular, the minimal element~$\hat 0$  and all atoms of $P$ are
irreducible elements in~$P$. For $x\,{\in}\,P$, we call 
$D(x):= {\rm max}\,(I(P)_{\leq x})$ 
the {\em set of elementary divisors\/} of~$x$ -- a~term which is explained 
by the following proposition: 

\begin{prop} \label{prop_elem_div} 
Let $P$ be a poset with a unique minimal element~$\hat 0$.
For $x\,{\in}\,P$ there exists a unique finest decomposition of the 
interval~$[\hat 0,x]$ in~$P$ as a direct product, which is given 
by an isomorphism 
$\varphi_x^{\rm el}:\,  \prod_{j=1}^l\, [\hat 0,y_j] \,  
                         \stackrel{\cong}{\lra} \,  [\hat 0,x]
$,
with $\varphi_x^{\rm el}(\hat 0, \ldots, y_j, \ldots, \hat 0)\,{=}\, y_j$ 
for $j=1,\ldots,l$. The factors of this decomposition are the
intervals below the elementary divisors of $x$: $\{y_1,\ldots,y_l\}=D(x)$. 
\end{prop}

\pr 
Whenever a poset with a minimal element~$\hat 0$ is represented as a~direct 
product, all elements which have more than one coordinate different from 
$\hat 0$ are reducible. 
Hence, if $\prod_{j=1}^l[\hat 0,y_j]\,{\cong}\,[\hat 0,x]$, 
and the $y_j$ are irreducible for $j\,{=}\,1,\ldots,l$, then 
$\{y_1,\ldots,y_l\}\,{=}\, D(x)$.
\qed


\subsection{Building Sets}\label{ssect_buildg}
\mbox{ } 

\noindent 
 In this subsection we define the notion of building sets of
 a~semilattice and develop their structure theory.

\begin{df} \label{df_buildg}
Let $\cl$ be a semilattice. A subset $\cg$ in~$\cl\,{\setminus}\,\{\hat 0\}$ 
is called 
a {\bf building set} of~$\cl$ if for any 
$x\,{\in}\,\cl\,{\setminus}\,\{\hat 0\}$  and
{\rm max}$\, \cg_{\leq x}=\{x_1,\ldots,x_k\}$ there is an isomorphism
of posets
\begin{equation}\label{eq_buildg}
\varphi_x:\,\,\, \prod_{j=1}^k\,\,\, [\hat 0,x_j] \,\, 
                               \stackrel{\cong}{\lra}
                                \,\, [\hat 0,x]
\end{equation}
with $\varphi_x(\hat 0, \ldots, x_j, \ldots, \hat 0)\, = \, x_j$ 
for $j=1,\ldots, k$. 
We call $F(x):=\max \cg_{\leq x}$  the {\bf set of factors} of $x$ in $\cg$.
\end{df}


 The next proposition provides several equivalent conditions for a subset
 of $\cl\,{\setminus}\,\{\hat 0\}$ to be a~building set.

\begin{prop} \label{prop_buildingsets}
For a semilattice $\cl$ and a subset $\cg$ of~$\cl\,{\setminus}\,\{\hat 0\}$ 
the following are equivalent:
\begin{itemize}

\item[(1)] $\cg$ is a building set of $\cl$;

\item[(2)] $\cg\/{\supseteq}\,I(\cl)\,{\setminus}\,\{\hat 0\}$, 
and for every $x\,{\in}\,\cl\,{\setminus}\,\{\hat 0\}$ with
$D(x)\,{=}\,\{y_1,\ldots, y_l\}$ the elementary divisors of~$x$, there 
exists a 
partition $\pi_x\,{=}\,\pi_1|\ldots|\pi_k$~of~$[l]$ 
with blocks $\pi_t\,{=}\,\{i_1,\ldots, i_{|\pi_t|}\}$ for 
$t\,{=}\,1,\ldots,k$, such that the elements 
in ${\rm max}\, \cg_{\leq x}\,{=} \,\{x_1,\ldots, x_k\}$ are of the form 
$
x_t\, = \, \varphi_x^{{\rm el}}(\hat 0, \ldots, \hat 0,
                        y_{i_1},\hat 0, \ldots, \hat 0,$ $
                        y_{i_2},\hat 0, \ldots, \hat 0, 
                        y_{i_{|\pi_t|}},\hat 0)
$.

\noindent
Informally speaking, the factors of $x$ in $\cg$ are products of disjoint
sets of elementary divisors of~$x$.
\item[(3)] $\cg$ generates $\cl\,{\setminus}\,\{\hat 0\}$ by $\vee$, and
for any $x\,{\in}\,\cl$, any $\{y,y_1,\dots,y_t\}\subseteq 
\max \cg_{\leq x}$, and  $z\,{\in}\,\cl$ with $z\,{<}\,y$, we have 
$
  \cg_{\leq y}\cap \cg_{\leq z\vee y_1\vee\dots\vee y_t}= \cg_{\leq z}
$.

\item[(4)] $\cg$ generates $\cl\,{\setminus}\,\{\hat 0\}$ by $\vee$, and for any $x\,{\in}\,\cl$, 
any $\{y,y_1,\dots,y_t\}\subseteq\max \cg_{\leq x}$,  and $z\,{\in}\,\cl$ with 
$z\,{<}\,y$,  the following two conditions are satisfied:
\[ 
\begin{array}{rll}
i)  & \cg_{\leq y}\cap \cg_{\leq y_1\vee\dots\vee y_t}=\emptyset
&\text{``disjointness,''}\\
i\!i) &z\vee y_1\vee\dots\vee y_t<y\vee y_1\vee\dots\vee y_t
&\text{``necessity.''} 
\end{array} 
\]
\end{itemize}
\end{prop}

\pr 

\noindent
\underline{(1)$\Rightarrow$(2)}: That $\cg$ contains 
$I(\cl)\,{\setminus}\,\{\hat 0\}$ follows directly
from the definition of building sets. We have the following isomorphisms:
$\varphi_x:\,\prod_{j=1}^k\, [\hat 0,x_j]\,{\lra}\,[\hat 0,x]$
by the building set property, and
$\varphi_{x_j}^{\rm el}:\,\prod_{y\in D(x_j)}\, [\hat 0,y]
   \,{\lra}\,
[\hat 0,x_j]$
for $j\,{=}\, 1,\ldots, k$ by Proposition~\ref{prop_elem_div}.
The composition 
$\varphi_x\,\circ\, (\prod_{j=1}^k\,\varphi_{x_j}^{\rm el})$
yields the finest decomposition $\varphi_{x}^{\rm el}$ of~$[\hat 0,x]$.
Thus, $D(x)\,{=}\,\uplus_{j=1}^k\, D(x_j)$, which gives the partition 
described in~(2).


\noindent
\underline{(2)$\Rightarrow$(1)}: The decomposition of~$[\hat 0,x]$ into
intervals below the elements in $\max \cg_{\leq x}$ follows from 
Proposition~\ref{prop_elem_div} by
assembling factors~$[\hat 0, y_j]$ with maximal elements indexed by elements 
from the same block of the partition~$\pi_x$ into one factor.

\noindent
\underline{(1)$\Rightarrow$(3)}: (3) is a direct consequence of~$[\hat 0,x]$
decomposing into a direct product of the form described in the definition of
building sets. 

\noindent
\underline{(3)$\Rightarrow$(4)}: $i)$ follows by setting $z=\hat 0$ 
in~(3).
Equality in $i\!i)$ implies with~(3) that $\cg_{\leq y}=\cg_{\leq z}$,
in particular, $y\in \cg_{\leq z}$ -- a contradiction to $z<y$.

\noindent
\underline{(4)$\Rightarrow$(1)}: For $x\in \cl\,{\setminus}\,\{\hat 0\}$ and max$\,\cg_{\leq x}\, = \, 
\{x_1,\ldots, x_k\}$ consider the poset map
\[
\phi:\,\,\, \prod_{j=1}^k\,\,\, [\hat 0,x_j] \,\, \lra \,\, [\hat 0,x] \, ,
\quad
      (\alpha_1, \ldots, \alpha_k)\,\,  \lmt \,\,  
                               \alpha_1 \vee \ldots  \vee \alpha_k\, .
\] 

\noindent
i) $\phi$ {\em is surjective\/}: For $\hat 0\,{\neq}\, y \,{\leq}\, x$, let 
max$\,\cg_{\leq y}=\{y_1,\ldots,y_t\}$. First, $\bigvee_{i=1}^t y_i\,{=}\,y$,
since $\cg$ generates $\cl$ by $\vee$. Second, define 
$\gamma_j\,{:=}\, \bigvee_{y_i\in S_j}y_i$ 
with 
$S_j\,{:=}\, ({\rm max}\,\cg_{\leq y})\,{\cap}\,\cg_{\leq x_j}$
for $j\,{=}\,1,\ldots,k$. 
Clearly, $\gamma_j\,{\in}\,[\hat 0,x_j]$, and 
$\cup_{j=1}^k S_j\,{=}\, {\rm max}\,\cg_{\leq y}$, 
since $\cg_{\leq y}\,{\subseteq}\,
\cg_{\leq x}$. Hence, $\phi(\gamma_1,\ldots, \gamma_k)=\bigvee_{i=1}^t y_i=y$.


\noindent 
ii) $\phi$ {\em is injective\/}: 
a) Assume $\phi(\alpha_1,\ldots,\alpha_k)\, {=}\, 
\phi(\beta_1,\ldots,\beta_k)= y \neq x$. Let 
max$\,\cg_{\leq y}=\{y_1,\ldots,y_t\}$. By induction on the number of elements
in~$[\hat 0,x]$ we can assume that $[\hat 0,y]$ decomposes as a direct product
$
[\hat 0,y] \, \cong \,  \prod_{i=1}^t\, [\hat 0,y_i] 
$.
Moreover, the subsets $S_j$ of max$\,\cg_{\leq y}$ defined in~i) actually
partition  max$\,\cg_{\leq y}$ as follows from the disjointness property
applied to pairwise intersections of the  $\,\cg_{\leq x_j}$.
Thus,
$
[\hat 0,y] \, \cong \,  \prod_{j=1}^k\, [\hat 0,\gamma_j] 
$,
with elements $\gamma_j\,{\in}\, [\hat 0,x_j]$ as above, and it follows that
$\alpha_j=\beta_j=\gamma_j$ for $j=1,\ldots,k$. \newline
b) Assume that $\phi(\alpha_1,\ldots,\alpha_k)\, {=}\, 
\phi(\beta_1,\ldots,\beta_k)\,{=}\,x$. By the necessity property it follows 
that $\alpha_j=\beta_j=x_j$ for $j=1,\ldots,k$. 
\qed

\vskip3pt
\begin{remrm} \label{rem_indpce_of_joins}
The definition of building sets and of irreducible elements, 
as well as the characterization of building sets
in Proposition~\ref{prop_buildingsets}~(2), are independent of the existence 
of a join operation and can be formulated for any poset with a unique minimal 
element. 
\end{remrm}

\noindent
We gather a few important properties 
of building sets. 

\pagebreak[4]
\begin{prop}\label{prop_buildingsets_properties}
For a building set $\cg$ of $\cl$, the following holds:
\begin{itemize}
\item [(1)] Let $x\in\cl$, $F(x)=\{x_1,\dots,x_k\}$ the set of factors of~$x$ 
in~$\cg$, and $\hat 0\,{\neq}\,y\,{\in}\,\cg$ with $y\,{\leq}\,x$. Then 
there exists a unique  $j\,{\in}\, \{1,\ldots,k\}$ such that 
$y\,{\leq}\,x_j$; i.e., $F(x)={\rm max}\,\cg_{\leq x}$ induces a partition of 
$\cg_{\leq x}$.
\item[(2)] For $x\,{\in}\,\cl$ and $x_0\,{\in}\,F(x)$,
\[
\bigvee\, (F(x)\,{\setminus}\,\{x_0\}) \,\, < \,\, \bigvee\, F(x)\, =x\, ,
\]
i.e., each factor of~$x$ in $\cg$ is needed to generate~$x$.
\item [(3)] If $h_1,\dots,h_k$ in $\cg$ are such that 
$(h_i,\bigvee_{j=1}^k h_j] \cap \cg=\emptyset$ for $i=1,\dots,k$, 
then $F(\bigvee_{j=1}^k h_j)\, = \, \{h_1,\ldots, h_k\}$.
\end{itemize}
\end{prop}

\pr
(1) is a consequence of Proposition~\ref{prop_buildingsets}~(4){\em i}), 
as was noted already in the proof of (4)$\Rightarrow$(1), part ii)~a),
in the previous proposition. 
Taking the full set of factors and setting~$z\,{=}\,\hat 0$ in 
Proposition~\ref{prop_buildingsets}~(4){\em i$\!$i}), yields~(2). For~(3) note
that $\{h_1,\ldots,h_k\}\,{\subseteq}\,F(\bigvee_{j=1}^k h_j)$ by assumption.
If $\{h_1,\ldots, h_k\}$ were not the complete set of factors, 
we would obtain a contradiction to~(2).
\qed

\begin{explrm} \label{expl_building}
(1) For the boolean lattice $\cb_n$ of rank~$n$, its atoms form the 
minimal building set. As with any other semilattice,
the full poset without its minimal element gives the maximal 
building set. 

In the smallest interesting example, the rank~$3$ boolean lattice $\cb_3$,
 we see that there are 
other building sets between these extremal choices: The atoms can be 
combined with any other rank~$2$ element to form a building set.
Moreover, atoms can be combined with the top element to form a 
building set,  and any other subset of $\cb_3$ containing the latter 
is in fact a building set.

\medskip
\noindent 
(2) For the partition lattice $\Pi_n$, the minimal building set 
is given by the \mbox{$1$-block} partitions. 
Again, the maximal building set 
is given by the full lattice without its minimal element. 
Looking at $\Pi_4$, we see that we can add any
$2$-block partition to the minimal building set, e.g., (12)(34), 
to obtain building sets other than the extreme ones.

\medskip
\noindent 
(3) The lattice $D_n$ of positive integral divisors of a natural
number~$n>0$ ordered by division relation has the prime powers
dividing $n$ as its minimal building set. Note that
this example includes the boolean lattice for any $n$ having no square
divisors, hence there are ample building sets between the extreme
choices.
\end{explrm}


\subsection{Nested Sets}
\mbox{ } 

\noindent 
  In this subsection we define the notion of nested subsets of
a~building set of a~semilattice and prove some of their properties.

\begin{df} \label{df_nested}
Let $\cl$ be a semilattice and $\cg$ a building set of $\cl$.
A subset $N$ in $\cg$ is called~{\bf nested} if, for any 
set of incomparable elements 
$x_1,\dots,x_t$ in $ N$ of cardinality at least two, 
the join $x_1\vee\dots\vee x_t$ exists and does not belong to $\cg$.
The nested sets in $\cg$ form an abstract simplicial complex, denoted 
$\cn(\cg)$. 
\end{df}

Note that the elements of $\cg$ are the vertices of the complex 
of nested sets~$\cn(\cg)$. Moreover, 
the order complex of $\cg$ is a subcomplex of $\cn(\cg)$, since
linearly ordered subsets of~$\cg$ are nested.

\begin{prop}\label{prop_nested}
For a given semilattice $\cl$ and a subset $N$ of a building set $\cg$ 
of~$\cl$, the following are equivalent:
  \begin{enumerate}
  \item [(1)] $N$ is nested.
  \item [(2)] Whenever $x_1,\dots,x_t$ are noncomparable elements in~$N$, 
the join $x_1\vee\dots\vee x_t$ exists, and 
$F(x_1\vee\dots\vee x_t)=\{x_1,\dots,x_t\}$.
  \item [(3)] There exists a chain $C\subseteq\cl$, such that
$N=\bigcup_{x\in C}F(x)$. 
  \item [(4)] $N\in\Lambda$, where $\Lambda$ is the maximal subset of
  $2^{\cg}$, for which the following three conditions are satisfied:
    \begin{enumerate}
    \item [(o)] $\emptyset\in\Lambda$, and $\{g\}\in\Lambda$, for $g\in \cg$;
    \item [(i)] if $N\in\Lambda$ and $x\in\max N$, then $N_{<x}\in\Lambda$;
    \item [(ii)] if $N\in\Lambda$, then $\max N=F(\bigvee \, \max N)$.
    \end{enumerate}
  \end{enumerate}
\end{prop}
\pr 

\noindent
\underline{(1)$\Rightarrow$(2)}:
Let $N$ be a nested set, and 
$M\,{=}\,\{x_1,\dots,x_t\}\,{\subseteq}\,N$  
a set of incomparable elements with $\bigvee_{i=1}^t x_i\,{\not \in}\,\cg$. 
We can assume that for some $x_j$: 
$(x_j,\bigvee_{i=1}^t x_i]\,{\cap}\,\cg\,{\neq}\,\emptyset$, 
otherwise the claim follows by 
Proposition~\ref{prop_buildingsets_properties}~(3).
Without loss of generality, we assume that there exists an element
$y\,{\in}\, (x_1,\bigvee_{i=1}^t x_i]\,{\cap}\,\cg$ and that 
$y\,{\in}\,\max \cg_{\leq \bigvee M}$. 
Define $M'\,{:=}\,\{x_1,\dots,x_t\}\,{\cap}\,\cg_{\leq y}\,{=}\,
\{x_1\,{=}\,x_{j_0},x_{j_1},\dots,x_{j_k}\}$ and 
$z\,{:=}\,\bigvee_{l=0}^k x_{j_l}$. 
Since $M'\,{=}\,\{x_{j_0},x_{j_1},\dots,x_{j_k}\}$ is nested
(it is a subset of $N$), we have the strict inequality $z\,{<}\,y$.
Furthermore,
\[
\bigvee_{i=1}^t x_i  \, = \, 
z\vee\bigvee (M\setminus M')\, \leq \, 
z\vee\bigvee(\max G_{\leq\bigvee M}\setminus\{y\})  \, < \, 
\bigvee_{i=1}^t x_i\,,
\]
where the first inequality follows from 
Proposition~\ref{prop_buildingsets_properties}~(1) and the second 
inequality from Proposition~\ref{prop_buildingsets_properties}~(2). 
We thus arrive to a contradiction, which finishes the proof.

\noindent
\underline{(2)$\Rightarrow$(1)}: Obvious.

\noindent
\underline{(2)$\Rightarrow$(3)}: Let $N$ be a~set satisfying condition~(2). 
Fix a~particular linear extension $\{x_1,\dots,x_k\}$ on the partial order 
of~$N$, and define
$\alpha_j\,{:=}\,x_1\,{\vee}\,\dots\,{\vee}\, x_j$, for $j\,{=}\,1,\dots,k$.
By (2) we have $F(\alpha_j)\,{=}\,\max\{x_1,\dots,x_j\}$, and therefore
$x_j\,{\in}\, F(\alpha_j)$ and $x_{j+1}\,{\not\in}\, F(\alpha_j)$ for
$j\,{=}\,1,\dots,k$. Hence, the $\alpha_j$'s are different and form
a~chain $C=\alpha_1<\alpha_2<\dots<\alpha_k$. By construction, 
$N=\bigcup_{x\in C}F(x)$.

\noindent
\underline{(1),(2)$\Rightarrow$(4)}: 
Let $N$ be a nested set, we shall prove that $N\in\Lambda$ 
by induction on the size of $N$:
\begin{enumerate}
\item if $|N|=0$, then $N\in\Lambda$ by condition (o);
\item if $|N|\geq 1$, then $\max N=F(\bigvee\,\max N)$ by 
condition~(2). Furthermore, since $|N_{<x}|<|N|$, and $N_{<x}$ 
is nested (it is a subset of $N$), $N_{<x}\in\Lambda$ by induction. 
Hence $N\in\Lambda$.
\end{enumerate}

\noindent
\underline{(3)$\Rightarrow$(1)}: 
Let $C\,{=}\,(\alpha_1\,{<}\,\dots\,{<}\,\alpha_k)$ be a chain in~$\cl$ and 
$N\,{=}\,\bigcup_{x\in C}F(x)$. 
Let $N'\,{=}\,\{x_1,\ldots,x_t\}\,{\subseteq}\, N$, $t\,{\geq}\,2$, be an
antichain in~$N$, and~$s$ the maximal index in~$C$ such that 
$N'\,{\cap}\,F(\alpha_s)\,{\neq}\,\emptyset$. In particular, 
$N'\,{\cap}\,F(\alpha_s)\,{\neq}\,\{\alpha_s\}$ due to $|N'|>1$ and
$N'$ being an antichain.

Let $y\,{\in}\,N'\,{\cap}\,F(\alpha_s)$. 
If $|N'\,{\cap}\,F(\alpha_s)|\,{>}\,1$,
\[
     y\, < \, \bigvee (N'\,{\cap}\,F(\alpha_s)) \, \leq \, 
              \bigvee \, N'      \, \leq \, \alpha_s\, ,  
\]
where the strict inequality is a consequence of the necessity property
for building sets. Thus, $ \bigvee  N'\,{\not \in}\, \cg$. 
If $|N'\,{\cap}\,F(\alpha_s)|\,{=}\,1$, we have 
$y\,{<}\,\bigvee N'\,{\leq}\,\alpha_s$,
due to $N'$ being an antichain with $|N'|\,{>}\,1$, and
again $\bigvee  N'\,{\not \in}\, \cg$.  

%

\noindent
\underline{(4)$\Rightarrow$(3)}: 
We need the following fact:

\vskip3pt
{\bf Fact.} {\it If there are elements $x_1,\dots,x_t$ and $y_1,\dots,y_k$ 
in~$\cl$, such that $x_t\,{>}\,y_j$ for $j\,{=}\,1,\dots,k$, and 
$F(\bigvee_{i=1}^t x_i)\,{=}\,\{x_1,\dots,x_t\}$, and  
$F(\bigvee_{j=1}^k y_j)\,{=}\,\{y_1,\dots,y_k\}$, 
then 
$F(x_1\vee\dots\vee x_{t-1}\vee y_1\vee\dots\vee y_k)\,{=}\,
\{x_1,\dots,x_{t-1},y_1,\dots,y_k\}.$}

\vskip3pt
   Once the fact above is proved, one can derive (3) as follows:
For $N\in\Lambda$ we shall form a chain 
$C\,{=}\,(\alpha_1\,{<}\,\dots\,{<}\,\alpha_{|N|})$ 
such that $N\,{=}\,\bigcup_{i=1}^{|N|} F(\alpha_i)$.
Choose a linear extension $\{x_1,\dots,x_t\}$ of $N$.
%
%
Set $\alpha_t\,{=}\,\bigvee \, \max N$,
$\alpha_{t-1}\,{=}\,\bigvee \, \max (N\,{\setminus}\,\{x_t\})$,
$\alpha_{t-2}\,{=}\,\bigvee \, \max (N\,{\setminus}\,\{x_t, x_{t-1}\})$,
and so on. By (4)(ii), $F(\alpha_t)\,{=}\,\max N$. Applying (4)(i) to
$x_t\,{\in}\,\max N$, and (4)(ii) to $N_{<x_t}$, we obtain 
$F(\bigvee\, \max N_{<x_t})\,{=}\,\max N_{<x_t}$. With the fact above,
we conclude that $F(\alpha_{t-1})\,{=}\,\max (N\,{\setminus}\,\{x_t\})$,
and, using the same argument iteratively, we arrive to 
$N\,{=}\,\bigcup_{i=1}^t F(\alpha_i)$. 
%
%

\vskip3pt
\noindent
{\bf Proof of the fact.}
Set $\alpha\,{:=}\,x_1\,{\vee}\,\dots\,{\vee}\, x_{t-1}\,{\vee}\, 
                   y_1\,{\vee}\,\dots\,{\vee}\, y_k$.
Since~$\alpha\,{\leq}\,\bigvee_{i=1}^t\,x_i$, the factors of~$\alpha$
can be partitioned into groups of elements below the $x_i$ for $i\,{=}\,
1,\ldots,t$, by Proposition~\ref{prop_buildingsets_properties}~(1).
Since $x_i\,{\leq}\,\alpha$ for $i\,{=}\,1,\ldots,t{-}1$, we obtain
$F(\alpha)\,{=}\,\{x_1,\ldots, x_{t-1}, \gamma_1,\ldots, \gamma_{m}\}$
with $\gamma_j\,{\leq}\,x_t$ for $j=1,\ldots,m$.

Again using  Proposition~\ref{prop_buildingsets_properties}~(1),
the $y_1,\dots, y_k$ can be partitioned into groups below the factors 
$\gamma_j$ for $j\,{=}\,1,\ldots,m$. The occurrence of one strict inequality
$\bigvee\,\{y_l\,|\,y_l\,{\leq}\,\gamma_j\}\,{<} \, \gamma_j$ 
for some $j\,{\in}\,
\{1,\ldots,m\}$ yields a contradiction to 
$\alpha\,{=}\, \bigvee_{i=1}^{t-1} x_i \,{\vee}$ $\bigvee_{j=1}^{k} y_j$
${=}\, \bigvee_{i=1}^{t-1} x_i \,{\vee}\,\bigvee_{j=1}^{m} \gamma_j
$, due to the necessity property of building sets. Moreover, since the
$y_i$ are factors themselves, joins of more than two of the $y_i$'s are
not elements of~$\cg$. Thus, $y_i\,{=}\,\gamma_i$, for 
$i\,{=}\,1,\ldots,k{=}m$, as claimed. \qed

\begin{explrm} \label{expl_nestedsets}
(1) For the boolean lattice $B_n$ with its minimal building set, 
any subset of atoms is nested.
The nested set complex hence is a simplex on $n$ vertices. As for 
any other semilattice with maximal building set, the nested sets 
are the totally ordered subsets of the poset, hence the nested set 
complex is the order complex of the poset. In the particular case of 
$\cb_n$ it is the barycentric subdivision of a simplex on~$n$ vertices.
For $\cb_3$ with building set $\cg=\{1,2,3,23\}$ the nested set complex
consists of two triangles, namely $\{1,2,23\}$ and $\{1,3,23\}$.

(2) For the partition lattice $\Pi_n$ with its minimal building set of
$1$-block partitions, a subset of such partitions is nested if and
only if any two non-trivial blocks are either contained one in another
or disjoint. This is the example which has suggested the terminology
of {\em nested\/} sets in the first place, it appeared as the central
combinatorial structure in the paper of Fulton~\& MacPherson~\cite{FM}
on models for configuration spaces of smooth complex varieties.
\end{explrm}


\section{Sequences of Combinatorial Blowups} 
\label{sect_comb_blowups}
$\,$
\noindent 
We introduce the notion of a combinatorial blowup of an~element
in a~semilattice and prove that the set of semilattices is closed
under this operation.

\subsection{Combinatorial Blowups}

\begin{df}
For a semilattice $\cl$ and an element $\al\in\cl$ we define 
a poset~$\bl_\al\cl$, the {\bf combinatorial blowup of $\cl$ at $\al$},  
as follows:
\begin{itemize}
\item[$\circ$] elements of $\bl_\al\cl$:
  \begin{enumerate}
  \item[(1)] $y\in\cl$, such that $y\not\geq\al$;
  \item[(2)] $[\al,y]$, for $y\in\cl$, such that $y\not\geq\al$
             and $(y\vee\al)_\cl$ exists \newline
             (in particular, $[\al,\hat 0]$ can be
             thought of as the result of blowing up $\al$);
  \end{enumerate} \pagebreak[4]
\item[$\circ$] order relations in $\bl_\al\cl$: 
  \begin{enumerate}
  \item[(1)] $y>z$ in $\bl_\al\cl$ if $y>z$ in $\cl$;
  \item[(2)] $[\al,y]>[\al,z]$ in $\bl_\al\cl$ if $y>z$ in $\cl$;
  \item[(3)] $[\al,y]>z$ in $\bl_\al\cl$ if $y\geq z$ in $\cl$;
  \end{enumerate}
  where in all three cases $y,z\not\geq\al$. 
\end{itemize}
\end{df}

Note that the atoms in $\bl_\al\cl$ are the atoms of~$\cl$
together with the element~$[\alpha,\hat 0]$.
It is easy, albeit tedious, to check that the class of 
(meet-)semilattices is closed under combinatorial blowups.

\begin{lm}
Let $\cl$ be a semilattice and $\al\in\cl$, then
$\bl_\al\cl$ is a semilattice.
\end{lm}
\pr The joins in $\bl_\al\cl$ are defined by the rule
\begin{align*}
([\al,y]\vee[\al,z])_{\bl_\al\cl}&=[\al,(y\vee z)_\cl],\\
  ([\al,y]\vee z)_{\bl_\al\cl}   &=[\al,(y\vee z)_\cl],\\
  (y\vee z)_{\bl_\al\cl}         &=(y\vee z)_\cl,
\end{align*}
which is applicable only if $(y\vee z)_\cl$ exists, otherwise
the corresponding joins in $\bl_\al\cl$ do not exist. Also,
the first and second formulae are applicable only in the case 
$(y\vee z)_\cl\not\geq\al$, otherwise the corresponding
joins do not exist. 
  The check of this is straightforward and is left to the reader.

\qed
 
Observe that it is possible that $(x\vee y)_\cl$ exists, while 
$(x\vee y)_{\bl_\al\cl}$ does not.


\subsection{Blowing Up Building Sets} \mbox{ }

\noindent
In this subsection we prove that if one combinatorially blows up 
a~building set of a~semilattice in any chosen linear extension order,
then one ends up with the face poset of the simplicial complex of nested
sets of this building set.
The following proposition provides the essential step for the proof. 

\begin{prop}\label{prop_singleblow}
  Let $\cl$ be a semilattice, $\cg$ a building set of~$\cl$, and 
$\al\,{\in}\,\max \cg$. Then, $\wti \cg=(\cg\sm\{\al\})\cup\{[\al,\hat 0]\}$
is a building set of $\bl_\al\cl$.
Furthermore, the nested subsets of $\wti \cg$ are precisely
the nested subsets of $\cg$ with $\al$ replaced by $[\al,\hat 0]$.
\end{prop}

\pr 
It is easy to see that $\wti \cg$ is a building set of $\bl_\al\cl$.
Indeed, given $x\in\cl\sm\cl_{\geq\al}$, \eqref{eq_buildg} is obvious for 
$x\in\bl_\al\cl$, and, if $(x\vee \al)_\cl$ exists, it follows for 
$[\al,x]\in\bl_\al\cl$ from the identity
\[
[\hat 0,[\al,x]]_{\bl_\al\cl}\, =\, [\hat 0,x]_{\bl_\al\cl}\times B_1\, ,
\]
where $B_1$ is the subposet consisting of the two comparable elements
$\hat 0\,{<}\,[\al,\hat 0]$.

Let us now see that the sets of nested subsets of $\cg$ and 
$\wti \cg$ are the same when replacing $\al$ by $[\al,\hat 0]$:
 
Let $N$ be a nested set in $\cg$, not containing $\al$. For incomparable
elements $x_1,\ldots,x_t$ in~$N$, 
$\bigvee_{i=1}^t\,x_i\,{\not \geq}\,\al$,
since otherwise we had 
$\al\,{\in}\, \max \cg_{\leq \bigvee x_i} =
F(\bigvee_{i=1}^tx_i)\,{=}\,
\{x_1,\ldots,x_t\}$
by Proposition~\ref{prop_nested}(2). Thus, $\bigvee_{i=1}^t x_i$ exists in
$\bl_{\al} \cl$ and $\bigvee_{i=1}^t x_i \,{\not\in}\,  
\wti \cg$. Hence, $N$ is nested in $\wti \cg$.
A nested subset in  $\wti \cg$ not containing $[\al,\hat 0]$ is obviously
nested in $\cg$.

Let now $N$ be nested in~$\cg$ containing $\alpha$, and set 
$\wti N\,{=}\,(N\sm\{\al\})\cup\{[\al,\hat 0]\}$. Subsets of incomparable 
elements in $\wti N$ not containing $[\al,\hat 0]$ can be dealt with as above.
Thus assume that $[\al,\hat 0],x_1,\ldots,x_t$ are incomparable in $\wti N$.
Then, $x_1,\ldots,x_t$ are incomparable in the nested set~$N$, and, as above,
we conclude that $\bigvee_{i=1}^t\,x_i$ exists and 
$\bigvee_{i=1}^t\,x_i\,{\not \geq}\,\al$. Moreover, $\alpha \vee 
\bigvee_{i=1}^t\,x_i$ exists in $\cl$ (joins of nested sets always exist!),
thus, 
$[\al,\bigvee_{i=1}^t\,x_i]\,{=}\,[\al,\hat 0]\vee \bigvee_{i=1}^t\,x_i$
exists in $\bl_{\al}\cl$ and is obviously not contained in $\wti \cg$.
We conclude that $\wti N$ is nested in $\wti \cg$. 

Vice versa, let $\wti N$ be nested in~$\wti \cg$ containing $[\al,\hat 0]$,
and set $N\,{=}\,(\wti N\sm \{[\al,\hat 0]\})\cup\{\al\}$.
Again it suffices to consider subsets of incomparable elements
$\al,x_1,\ldots,x_t$ in $N$. With  $[\al,\hat 0],x_1,\ldots,x_t$ 
incomparable in $\wti N$, $[\al,\hat 0]\vee \bigvee_{i=1}^t\,x_i\,{=}\,
[\al,\bigvee_{i=1}^t\,x_i]$ exists in $\bl_{\al}\cl$, thus 
$\al \vee \bigvee_{i=1}^t\,x_i$ exists in~$\cl$. Incomparability
implies that $\al \vee \bigvee_{i=1}^t\,x_i>\al$, and thus
$\al \vee \bigvee_{i=1}^t\,x_i \not \in \cg$. We conclude that~$N$ 
is nested in~$\cg$.
\qed




\medskip
By iterating the combinatorial blowup described in 
Proposition~\ref{prop_singleblow} through all of~$\cg$, 
we obtain the following theorem, which serves as a motivation 
for the entire development.

\begin{thm} \label{thm_main}
Let $\cl$ be a semilattice and $\cg$ a building set of~$\cl$
with some chosen linear extension: $\cg=\{G_1,\dots,G_t\}$, with  
$G_i>G_j$ implying $i<j$. Let $\bl_k\cl$ denote the result 
of subsequent blowups $\bl_{G_k}(\bl_{G_{k-1}}(\dots\bl_{G_1}\cl))$.
Then the final semilattice $\bl_t\cl$ is equal to the face poset of 
the simplicial complex $\cn(\cg)$.
\end{thm}

\pr
The building set~$\cg_t$ of $\bl_t\cl$ that results from iterated 
application of Proposition~\ref{prop_singleblow}  obviously is the set
of atoms $\mathfrak A$ in $\bl_t\cl$. Every element $x\in\bl_t\cl$ is the join 
of atoms below it: $x=\bigvee \mathfrak A_{\leq x}$. The subset 
$\mathfrak A_{\leq x}$
of $\cg_t$ is nested, in particular, it is the set of factors of~$x$ 
in $\bl_t\cl$ with respect to $\cg_t$ (Proposition~\ref{prop_nested}(2)). 
Proposition~\ref{prop_buildingsets_properties}(2) implies that the
interval~$[\hat 0,x]$ in $\bl_t\cl$ is boolean. We conclude that
$\bl_t\cl$ is the face poset of a simplicial complex with faces in 
one-to-one correspondence with the nested sets in $\cg_t$, which in turn
correspond to the nested sets in~$\cg$ by Proposition~\ref{prop_singleblow}.
\qed


\section{Instances of combinatorial blowups}
\label{sect4}

\subsection{De Concini-Procesi Models Of Subspace Arrangements}
\label{ssect_DPmodels}
$\,$
\vskip3pt

Let $\ca=\{A_1,\ldots,A_n\}$ be an arrangement of linear subspaces 
in complex space~$\dc^d$. Much effort has been spent on describing the 
cohomology of the complement~$\cm(\ca)=\dc^d\,{\setminus}\,\bigcup \ca$
of such an arrangement and, in particular, on answering the question 
whether the cohomology algebra is completely determined by 
the combinatorial data of the arrangement. Here, combinatorial data 
is understood as the lattice~$\cl(\ca)$ of intersections of subspaces 
of~$\ca$ ordered by reverse inclusion together with the complex 
codimensions of the intersections. A major step towards the solution 
of this problem (for a complete answer see~\cite{DGM, dLS}) was 
the construction of smooth models for the complement~$\cm(\ca)$ by 
De Concini~\& Procesi~\cite{DP1} that allowed for an explicit 
description of rational models for~$\cm(\ca)$ following~\cite{M}.  
The De Concini-Procesi models for arrangements in turn are one instance 
in a sequence of model constructions reaching from compactifications 
of symmetric spaces~\cite{DP3,DP4}, over the Fulton-MacPherson 
compactifications of configuration spaces~\cite{FM} to the general 
framework of wonderful conical compactifications proposed by 
MacPherson~\& Procesi~\cite{MP}.

Given a complex subspace arrangement~$\ca$ in~$\dc^d$, De Concini~\& Procesi
describe a smooth irreducible variety~$Y$ together with a proper map 
$\pi:\, Y\, \lra \, \dc^d$ such that $\pi$ is isomorphism over $\cm(\ca)$,
and the complement of the preimage of~ $\cm(\ca)$ is a union of irreducible 
divisors with normal crossings in~$Y$. The model~$Y$ can be
constructed by a~sequence of blowups of smooth subvarieties that is 
prescribed by the stratification of complex space induced by the arrangement. 


\subsubsection{Building sets for subspace arrangements}

In order to enumerate the strata in the intersection stratification of~$Y$
given by the irreducible divisors, De Concini \& Procesi introduced the 
notions of building sets, nested sets and irreducible elements as follows:  

\begin{df} \label{df_DPnotions}
{\rm (\cite[\S 2]{DP1})}
  Let $\cl(\ca)$ be the intersection lattice of an arrangement~$\ca$ 
of linear subspaces in a~finite dimensional complex vector space. 
Consider the lattice~$\cl(\ca)^*$ formed by the orthogonal 
complements of intersections ordered by inclusion. 
\begin{itemize}
\item[(1)] For $U\,{\in}\,\cl(\ca)^*$,
$U\, = \oplus_{i=1}^k\, U_i$ with $U_i\,{\in}\,\cl(\ca)^*$, is called a
{\bf decomposition\/} of~$U$ if for any $V\,{\subseteq}\,U$, 
$V\in \cl(\ca)^*$, $V\, = \oplus_{i=1}^k\, (U_i\cap V)$ and 
$U_i\cap V\,{\in}\,\cl(\ca)^*$ for $i=1,\ldots,k$. 
\item[(2)] Call  $U\,{\in}\,\cl(\ca)^*\,{\setminus}\,\{\hat 0\}$ 
{\bf irreducible\/} if it does not
admit a non-trivial decomposition.
\item[(3)] $\cg\,{\subseteq}\,\cl(\ca)^*\,{\setminus}\,\{\hat 0\}$ is called a 
{\bf building set\/}
for~$\ca$ if for any $U\,{\in}\,\cl(\ca)^*$ and $G_1,\ldots, G_k$
maximal in~$\cg$ below~$U$,
$U\, = \oplus_{i=1}^k\, G_i$ is a decomposition (the $\cg$-decomposition)
of~$U$.
\item[(4)] A subset $\cs\,{\subseteq}\,\cg$ is called {\bf nested\/} if for
any set of non-comparable elements $U_1,\ldots,U_k$ in~$\cs$, 
$U\, = \oplus_{i=1}^k\, U_i$ is the $\cg$-decomposition of~$U$.
\end{itemize}
\end{df}

Note that $\cl(\ca)^*$ coincides with $\cl(\ca)$ as abstract lattices.
We will therefore talk about irreducible elements, building sets and 
nested sets in $\cl(\ca)$ without explicitly referring to the dual setting
of the preceding definition.

The notions of Definition~\ref{df_DPnotions} are in part based on the 
earlier notions introduced by 
Fulton~\& MacPherson in~\cite{FM} to study compactifications 
of configuration spaces. Our terminology is naturally adopted 
from~\cite{FM,DP1}. Building sets and nested sets in the sense of 
De Concini~\& Procesi are building and nested sets for the intersection 
lattices of subspace arrangements in our combinatorial sense
(see Proposition~\ref{combgeomprop}~(1) below). However, there
are differences. The opposite is not true: A combinatorial building set 
for the intersection lattice of a subspace arrangement is not
necessarily a building set for this arrangement in the sense of 
De Concini~\& Procesi, neither
are irreducible elements in the sense of De Concini~\& Procesi 
irreducible in our sense.

\begin{explrm} \label{expl_combversusDP}
({\em Combinatorial versus De Concini-Procesi 
building sets}) \newline
Consider the following arrangement $\ca$ of 3 subspaces 
in~$\dc^4$:
\[
A_1:\,\, z_4  =  0\, , \quad
A_2:\,\, z_1= z_2 =  0\, ,\quad
A_3:\,\, z_1= z_3 =  0 \,.
\]
The intersection lattice $\cl(\ca)$ is a boolean algebra on 3
elements.  Combinatorial building sets of this lattice have been
discussed in Example~\ref{expl_building}, in particular, the set of
atoms $\{A_1,A_2,A_3\}\,{\subseteq}\,\cl(\ca)$ is the minimal
combinatorial building set. However, any building set for $\ca$ in the
sense of De Concini~\& Procesi necessarily includes the intersection
$A_2\,{\cap}\,A_3$, since its orthogonal complement does not decompose
in~$\cl(\ca)^*$. The minimal building set for~$\ca$, i.e., the set of
irreducibles for~$\ca$, in the sense of
De~Concini~\& Procesi is $\{A_1,A_2,A_3,A_2\,{\cap}\,A_3\}$. Any other
building set contains this minimal building set and 
the total intersection $\bigcap \ca\,{=}\,0$.
\end{explrm}

The main difference between our combinatorial set-up and the original 
context of De~Concini-Procesi model constructions can be formulated 
in the following way:
our constructions are order-theoretically canonical for a given 
semilattice. The set of combinatorial building sets, in particular 
the set of irreducible elements, depends only on the semilattice
itself and not on the geometry of the subspace arrangement which it 
encodes. See Proposition~\ref{combgeomprop} for a complete explanation.

\subsubsection{Local subspace arrangements}
 In order to trace the De Concini-Procesi construction step by step
we need the more general notion of a~local subspace arrangement.

\begin{df} \label{lsadf}
  Let $M$ be a smooth complex $d$-dimensional manifold and $\ca$ 
a~union of finitely many smooth complex submanifolds of $M$ such that 
all non-empty intersections of submanifolds in~$\ca$ are connected 
smooth complex submanifolds. $\ca$ is called a~{\bf local subspace 
arrangement}~if for any $x\in\ca$ there exists an~open set $N$ in~$M$ 
with $x\,{\in}\,N$, a~subspace arrangement $\wti\ca$ in ${\dc}^d$, 
and a~biholomorphic map $\phi:N\ra{\dc}^d$, such that 
$\phi(N\cap\ca)=\wti\ca$.
\end{df}

  Given a subspace arrangement $\ca$, the initial ambient space~$\dc^d$ 
of $\cm(\ca)$ carries a natural stratification  by the subspaces
of~$\ca$ and their intersections, the poset of strata being the 
intersection lattice $\cl(\ca)$ of the arrangement. For a~local
subspace arrangement $\ca=\{A_1,\dots,A_n\}$ in $M$ we again
consider the stratification of $M$ by all possible intersections
of the $A_i$'s, just like in the global case. The poset of strata
is also denoted by $\cl(\ca)$ and is called the intersection 
semilattice (it is a lattice if the intersection of all maximal
strata is nonempty).

\begin{df} \label{lsabsdf}
  Let $\ca$ be a local subspace arrangement and $\cl(\ca)$ its 
intersection semilattice. For $U\in\cl(\ca)$, $U_1,\dots,U_k\in\cl(\ca)$
are said to form a~{\bf decomposition} of $U$ if for any $x\in U$ there
exists an~open set $N$ with $x\,{\in}\,N$ and a~biholomorphic map 
$\phi:N\ra{\dc}^d$, such that $\phi(N\cap U_1),\dots,\phi(N\cap U_k)$
form a~decomposition of $\phi(N\cap U)$ in the sense of 
{\rm Definition~\ref{df_DPnotions}(1)}. 

  As in the global case, $\cg\subseteq\cl(\ca)$ is a~{\bf building 
set} for $\ca$ if for any $U\in\cl(\ca)$, the set of strata $\max\cg_{\leq U}$
gives a~decomposition of $U$.
\end{df}
 We shall refer to these building sets as {\em geometric\/} building sets.
The difference between combinatorial building sets and geometric
ones is contained in the dimension function as is explained in the 
following proposition.

\begin{prop} \label{combgeomprop} $\,$
  Let $\ca$ be a~local subspace arrangement with intersection 
semilattice~$\cl(\ca)$.
\begin{enumerate}
\item[(1)] If $\cg\subseteq\cl(\ca)$ is a~geometric building set of
$\ca$, then it is a~combinatorial building set.

\item[(2)] If $\cg\subseteq\cl(\ca)$ is a~combinatorial building set 
of $\cl(\ca)$, and for any $x\in\cl(\ca)$ the sum of codimensions of its
factors is equal to the codimension of $x$, then $\cg$ is
a~geometric building set.
\end{enumerate}
\end{prop}

\pr 
   In both cases it is enough to consider the case when $\ca$
is a~subspace arrangement.

(1) Consider~$\cg$ as a subset of~$\cl(\ca)^*$, then, for 
$U\,{\in}\,\cg$, the isomorphism $\varphi_U$ requested in 
Definition~\ref{df_buildg} is given by taking direct sums:
$$\varphi_U:\,\,\, \prod_{j=1}^k\,\,\, [\hat 0,G_j] \,\, 
\stackrel{\oplus_{j=1}^k}{\lra} \,\, [\hat 0,U]\, ,$$
where $G_1,\ldots,G_k$ are maximal in $\cg$ below~$U$.

(2) For $U\in\cl(\ca)^*$, the set $\{U_1,\dots,U_k\}=\max\cg_{\leq U}$
gives a~decomposition of $U$ because:
\begin{enumerate}
\item [a)] By the definition of $\cl(\ca)^*$ and the definition of
combinatorial building sets, we have $U=\Span(U_1,\dots,U_k)$, and, 
since $\sum_{i=1}^k\dim U_i=\dim U$, we have $U=\bigoplus_{i=1}^k U_i$;
\item [b)] for any $V\subseteq U$, $\bigoplus_{i=1}^k(U_i\cap V)
\subseteq V=\Span(U_1\wedge V,\dots,U_k\wedge V)\subseteq
\bigoplus_{i=1}^k(U_i\cap V)$, where "$\wedge$" denotes the meet 
operation in $\cl(\ca)^*$, hence $V=\bigoplus_{i=1}^k(U_i\cap V)$.
\qed
\end{enumerate}


\subsubsection{Intersection stratification of local arrangements 
               after blowup}
  Let a space~$X$ be given with an intersection stratification 
induced by a~local subspace arrangement, and let $G$ be a stratum in~$X$.
In the blowup of~$X$ at~$G$, $\bl_{G}X$, we find the following maximal strata: 
\begin{itemize}
\item[$\circ$] maximal strata in~$X$ that do not intersect with~$G$,
\item[$\circ$] blowups of maximal strata~$V$ at $G\cap V$, 
               $\bl_{G\cap V}V$,
               where $V$ is maximal in~$X$ and intersects~$G$,
\item[$\circ$] the exceptional divisor $\widetilde G$ replacing~$G$.
\end{itemize}
We consider the intersection stratification of $\bl_{G}X$ induced by
these maximal strata. We will later see (proof of 
Proposition~\ref{prop_compDP_cbl}) that in case $G$ is maximal
in a building set for the local arrangement in~$X$, then the union of
maximal strata in $\bl_GX$ is again a local arrangement with induced 
intersection stratification. In general, this is not the case,
see Example~\ref{expl_bllocarrgts}

For ease of notation, let us agree here that formally blowing up an
empty (non-existing) stratum has no effect on the space.
We think about a stratum~$Y$ in~$X$, intersection of all maximal
strata $V_1,\ldots,V_t$ that contain~$Y$, as being replaced by the
intersection of corresponding maximal strata in~$\bl_{G}X$:
\begin{equation} \label{eq_blownupstratum}
\bl_{G\cap V_1} V_1 \,\,  \cap \, \,
             \ldots \,\,  \cap \, \,          
\bl_{G\cap V_t} V_t\, ,
\end{equation}
(recall that $\bl_{G\cap V_j} V_j\,{=}\,V_j$ for
$G\,{\cap}\,V_j\,{=}\,\emptyset$). The
intersection~(\ref{eq_blownupstratum}) being empty means that the 
stratum~$Y$ vanishes under blowup of~$G$. For notational convenience,
we most often retain names of strata under blowups, thereby
referring to the replacement of strata described above. 

\begin{explrm} \label{expl_bllocarrgts}
({\em Local subspace arrangements are not closed under blowup})
\newline
We give an example which shows that blowing up a stratum in a local 
subspace arrangement does not necessarily result in a local subspace 
arrangement again. Consider the following arrangement of 2~planes and 
1~line in~$\dc^3$:
\[
A_1:\, \, y-z=0\, ,\quad
A_2:\, \, y+z=0\, ,\quad
L:\, \, x=y=0\, . 
\]
After blowing up~$L$, the planes $A_1$ and $A_2$ are replaced by 
complex line bundles over~$\dc {\rm P}^1$, which have in common their 
zero section~$Z$ and a complex line~$Y$; $L$ is replaced by a direct 
product of~$\dc$ and $\dc {\rm P}^1$, which intersects both line
bundles in~$Z$. The new maximal strata fail to form a local subspace 
arrangement in the point $Z\cap Y$.
\end{explrm}

\subsubsection{Tracing incidence structure during arrangement model 
construction} \label{ssect_tracg_inc}

  We now give a more detailed description of the model construction
by De Concini \& Procesi via successive blowups, and then proceed with 
linking our notion of combinatorial blowups to the context of arrangement 
models.

Let $\ca$ be a complex subspace arrangement,
$\cg\,{\subseteq}\,\cl(\ca)$ a~geometric building set for~$\ca$, and 
$\{G_1,\ldots,G_t\}$ some linear extension of the partial containment 
order on associated strata in~$\dc^d$ such that $G_k\,{\supset}\,G_l$ 
implies $l<k$. The De Concini-Procesi model $Y=Y_{\cg}$ of $\cm(\ca)$ 
is the result of blowing up the strata indexed by elements of~$\cg$ 
in the given order. Note that the linear order was chosen so that at 
each step the stratum which is to be blown up  does {\em not\/}
contain any other stratum indexed by an element of~$\cg$. At each step 
we consider intersection stratifications as described above, and we 
denote the poset of strata after blowup of $G_i$
with~$\cl^{\cg}_i(\ca)$. For the case of a~stratum $G_i$ being empty 
after previous blowups remember our agreement of considering blowups 
of $\emptyset$ as having no effect on a space. The later 
Proposition~\ref{prop_compDP_cbl} however shows that strata indexed 
by elements in~$\cg$ do not disappear during the sequence of blowups.

Let us remark that the combinatorial data of the initial
stratification, i.e., of the arrangement, prescribes much 
of the geometry of~$Y_{\cg}$: the complement $Y_{\cg}\,{\setminus}\,\cm(\ca)$
is a union of smooth irreducible divisors indexed by elements 
of~$\cg$, and these divisors intersect if and only if the set of
indices is nested in~$\cg$~\cite[Thm 3.2]{DP1}.


\begin{prop}\label{prop_compDP_cbl}
Let $\ca$ be an arrangement of complex subspaces, $\cg$ a building set
for $\ca$ in the sense of De Concini~\& Procesi, and $\{G_1,\ldots,G_t\}$ 
some linear extension of the partial containment order on associated strata
as described above. Let $\bl_i^{\cg}(\ca)$
denote the geometric result of successively blowing up strata 
$G_1,\ldots,G_i$, for $1\,{\leq}\,i\,{\leq}\,t$. Then,
\begin{itemize}
\item[(1)] 
The poset of strata $\cl^{\cg}_i(\ca)$  of $\bl_i^{\cg}(\ca)$ can be 
described as the result of a sequence of combinatorial blowups of the 
intersection lattice~$\cl\,{=}\,\cl(\ca)$:
\[
      \cl^{\cg}_i(\ca)\, \, = \, \, \bl_i(\cl)\, ,\qquad \mbox{for }\,
      1\,{\leq}\,i\,{\leq}\,t\, .
\]
{\rm (}Recall that $\bl_i(\cl)\,{=}\,\bl_{G_i}(\bl_{G_{i-1}}(\ldots \bl_{G_1}\cl))$
for $ 1\,{\leq}\,i\,{\leq}\,t$.{\rm )}
\item[(2)]
The union of maximal strata~$\ca_i^{\cg}$ in
$\bl_i^{\cg}(\ca)$ is a local subspace arrangement, with $\cg$ in 
$\cl^{\cg}_i(\ca)$ being a~building set for $\ca_i^{\cg}$ in the sense of \/ 
{\rm Definition~\ref{lsabsdf}}. {\rm (}Recall that $\cg$ here refers to 
the preimages 
of the original strata in $\cg\,{\subseteq}\,\cl(\ca)$ 
under the sequence of blowups.{\rm )}
\end{itemize}
\end{prop}

\noindent
\pr
We proceed by induction on the number of blowups. 
The induction start is obvious, since
the lattice of strata~$\cl_0^{\cg}(\ca)$ of the initial
stratification of $\dc^d$ coincides with the intersection lattice
$\cl(\ca)\,{=}\,\bl_0(\cl)$ of the arrangement~$\ca$. The union of maximal 
strata is the arrangement~$\ca$ itself with its given building set~$\cg$.

Assume that
$\cl_{i-1}^{\cg}(\ca)=\bl_{i-1}(\cl)$ for some $ 1\,{\leq}\,i\,{\leq}\,t$,
the union of maximal strata $\ca_{i-1}^{\cg}$ in $\bl_{i-1}^{\cg}(\ca)$ being
a local arrangement, and $\cg$ a building set for $\cl_{i-1}^{\cg}(\ca)$. 
Let $G\,{=}\,G_i$ be the next stratum to be blown up. 
First, we proceed in 4 steps to show that 
$\cl^{\cg}_i(\ca)\,{=}\,\bl_i(\cl)$. In 2 further steps we then verify
the claims in~(2).

\noindent
{\bf Step 1:} {\em Assign strata of $\bl_i^{\cg}(\ca)$ to elements 
in~$\bl_i(\cl)$.\/}\newline
We distinguish two types of elements in $\bl_i(\cl)$:
\[
\begin{array}{lrl}
\mbox{Type~I}:     & Y & 
                   \mbox{with }\, Y\,{\in}\,\bl_{i-1}(\cl) 
                   \mbox{ and}\, Y \not\geq G\, , \\
\mbox{Type~~I$\!$I}:\quad & [G,Y] & 
                   \mbox{with }\,Y\,{\in}\,\bl_{i-1}(\cl)\, , 
                   \, \, Y \not\geq G\, , \\
& &                \mbox{and }\, Y\,{\vee}\,G 
                   \mbox{ exists  in }\,\bl_{i-1}(\cl)\, . 
\end{array}
\]

\noindent
To $Y\,{\in}\,\bl_i(\cl)$ of type~I, assign $\bl_{G\cap Y}Y$ (recall
that blowing up an empty stratum does not change the space). Note that 
$\dim\bl_{G\cap Y}Y=\dim Y$.

\noindent
To $[G,Y]\,{\in}\,\bl_i(\cl)$ of type~I$\!$I, assign  
$(\bl_{G\cap Y}Y)\, \cap \, \widetilde G$, where $\widetilde G$ denotes
the exceptional divisor that replaces $G$ in $\bl_i^{\cg}(\ca)$. This 
description comprises $\widetilde G$ being assigned to $[G,\hat 0]$.
Note that $\dim(\bl_{G\cap Y}Y)\, \cap \, \widetilde G=\dim Y-1$.

\noindent
{\bf Step 2:} {\em Reverse inclusion order on the assigned spaces
                   coincides with the partial order
                   on~$\bl_i(\cl)$.\/}\newline
\noindent
(1) $X,Y\,{\in}\, \bl_i(\cl)$, both of type~I: 
\[
X\, \leq_{\bl_i(\cl)} Y \,\, \Leftrightarrow \,\,
X\, \leq_{\bl_{i-1}(\cl)} Y \,\, \Leftrightarrow \,\,
X \supseteq_{\bl_{i-1}^{\cg}(\ca)} Y \,\, \Leftrightarrow \,\,
\bl_{G\cap X} X \supseteq \bl_{G\cap Y} Y\, ,
\]
where ``$\Leftarrow$'' in the last equivalence can be seen by first
noting that $Y \,{\setminus}\, (G\,{\cap}\,Y) \subseteq X 
\,{\setminus}\, (G\,{\cap}\,X)$, and then comparing points in the
exceptional divisors.

\noindent
(2) $X, [G,Y]\,{\in}\, \bl_i(\cl)$, $X$ of type~I,
$[G,Y]$ of type~I$\!$I: \newline
As above we conclude
\begin{eqnarray*}
X\, \leq_{\bl_i(\cl)} [G,Y] & \Leftrightarrow & 
X\, \leq_{\bl_{i-1}(\cl)} Y  \\
                            & \Leftrightarrow &
X \supseteq_{\bl_{i-1}^{\cg}(\ca)} Y \,\, \Rightarrow \,\,
\bl_{G\cap X} X \supseteq \bl_{G\cap Y} Y \cap \widetilde G\, .
\end{eqnarray*}
To prove the converse is rather subtle. Note first that 
$G\,{\cap}\,Y \subseteq G\,{\cap}\,X$. 
Assume that $G$ strictly contains~$G\,{\cap}\,X$, then 
both $G\,{\cap}\,X$ and $G\,{\cap}\,Y$ are not in the building set due 
to the linear order chosen on~$\cg$, and~$G$ is a factor of both
$G\,{\cap}\,X$ and $G\,{\cap}\,Y$. Let 
$F(G\,{\cap}\,X)\,{=}\,\{G,G_1,\ldots,G_k\}$, 
$F(G\,{\cap}\,Y)\,{=}\,\{G,H_1,\ldots,H_t\}$. $X$ written as a join
of elements in~$\bl_{i-1}(\cl)$ below the factors of~$G\,{\cap}\,X$
reads
\[
X \, = \, g_X \vee  Z_1  \vee \ldots  \vee Z_k
\] 
for some $g_X\,{\in}\,[\hat 0,G]$, $Z_i\,{\in}\,[\hat 0, G_i]$ for 
$i\,{=}\, 1,\ldots, k$. If $Z_i\,{<}\,G_i$ for some $i{\in}\{1,\ldots,k\}$,
we have
\begin{eqnarray*}
G\vee X 
& = & 
G \vee (g_X \vee Z_1 \vee  \ldots \vee Z_i  \vee \ldots  \vee Z_k) \\ 
& \leq & 
G \vee (g_X \vee G_1 \vee  \ldots \vee Z_i  \vee \ldots  \vee G_k) \\
& < & 
G \vee  G_1 \vee  \ldots \vee G_k \, = \,  G \vee X\, ,  
\end{eqnarray*}
by the ``necessity'' property of 
Proposition~\ref{prop_buildingsets}(4), yielding a contradiction.
Hence,
\[
X \, = \, g_X \vee  G_1  \vee \ldots  \vee G_k\, ,
\]
and similarly, $Y\,{=}\, g_Y \vee H_1 \vee\ldots \vee H_t$
for some $g_Y\,{\in}\,[\hat 0,G]$.

For each $j\,{\in}\,\{1,\ldots,k\}$ there exists a unique 
$i_j\,{\in}\,\{1,\ldots,t\}$ such that $G_j\leq H_{i_j}$ by 
Proposition~\ref{prop_buildingsets_properties}(1). 
Thus, $\bigvee G_i\,{<}\, \bigvee H_j$, and, for showing that 
$X\leq Y$, it is enough to see that~$g_X\leq g_Y$.

We show that in an open neighborhood of any point
$y\,{\in}\,G\,{\cap}\,Y$, $g_Y\,{\subseteq}\,g_X$. This yields our
claim since strata in~$\bl_{i-1}^{\cg}(\ca)$ have pairwise transversal
intersections: if they coincide locally, they must coincide globally. 
With $\ca_{i-1}^{\cg}$ being a local arrangement, there exists
an open neighborhood of $y\,{\in}\,G\,{\cap}\,Y$  where
the stratification is biholomorphic to a stratification induced by a
subspace arrangement. We tacitly work in the arrangement setting,
using that $(\bl_{i-1}(\cl))_{\leq G\vee Y}$ is the intersection lattice 
of a product arrangement. The $\cg$-decomposition of $(G\vee Y)^{\perp}$
described in Definition~\ref{lsabsdf} yields (when transferred to the 
primal setting):
\[
g_Y\, = \,\Span(G,Y)\, .
\]
Analogously, $g_X\,{=}\,\Span(G,X)$.

In the linear setting we are concerned with, we interpret points in
the exceptional divisor of a blowup as follows:
\begin{equation}\label{eq_blowup}
\bl_{G\cap Y} Y \cap \widetilde G\,\, = \,\,
\{ (a,\Span(p,G\cap Y)) \, | \, 
a\in G\cap Y,\,  p\in Y\setminus (G\cap Y)  \}\, . 
\end{equation}
In terms of this description, the inclusion map 
$\bl_{G\cap Y} Y \cap \wti G \hookrightarrow \bl_G(\bl_{i-1}^{\cg}(\ca))$ reads
\[
  (a,\Span(p,G\cap Y)) \longmapsto  (a,\Span(p,G))\, .
\]
Therefore, $\bl_{G\cap Y} Y \cap \wti G$ being contained in
$\bl_{G\cap X} X\,{\subseteq}\, \bl_G(\bl_{i-1}^{\cg}(\ca))$ 
means that for $(a,\Span(p,G\cap Y))\in 
\bl_{G\cap Y} Y \,{\cap}\, \widetilde G$ there exists 
$q\,{\in}\, X\setminus (G\cap X)$ such that $\Span(p,G)\, = \,
\Span(q,G)$. In particular, $\Span(Y,G)\,{\subseteq}\,\Span(X,G)$, which
by our previous arguments implies that $Y\,{\subseteq}\,X$.

We assumed above that $G\,{\supset}\,G\,{\cap}\,X$. If $G\,{\cap}\,X$
coincides with~$G$, i.e., $X$ contains~$G$, then $g_X\,{=}\,X$ and a
similar reasoning applies to see that $Y\,{\subseteq}\,X$. Similarly
for $G\,{\cap}\,X\,{=}\,G\,{\cap}\,Y\,{=}\,G$.

\noindent
(3) $[G,X]$, $[G,Y]\,{\in}\, \bl_i(\cl)$, both of type~I$\!$I: 
\begin{eqnarray*}
[G,X]\, \leq_{\bl_i(\cl)} [G,Y] & \Leftrightarrow & 
X\, \leq_{\bl_{i-1}(\cl)} Y  \\
                            & \Leftrightarrow &
X \supseteq_{\bl_{i-1}^{\cg}(\ca)} Y \,\, \Leftrightarrow \,\,
\bl_{G\cap X} X \cap \widetilde G \supseteq \bl_{G\cap Y} Y \cap \widetilde G\, ,
\end{eqnarray*}
where ``$\Leftarrow$'' follows from~(2) and 
$\bl_{G\cap X} X\,{\supseteq}\,
\bl_{G\cap X} X \cap \widetilde G\,{\supseteq}\,
 \bl_{G\cap Y} Y \cap \widetilde G$.

\noindent
{\bf Step 3:} {\em Each of the assigned spaces is the intersection of 
maximal strata in~$\bl_i^{\cl}(\ca)$.}\newline
It is enough to show that spaces assigned to elements of type~I
in~$\bl_i(\cl)$ are intersections of new maximal strata. Those
associated to elements of type~I$\!$I then are intersections as well
by definition.

Let $Y\,{\in}\,\bl_i(\cl)$, $Y\,{\not\geq}\,G$, 
and $Y\,{=}\,\cap_{i=1}^t\, V_i$ with $V_1,\ldots,V_t$ the maximal strata 
in~$\bl_{i-1}^{\cg}(\ca)$ containing~$Y$. We claim that
\begin{equation} \label{eq_intdiv}
\bl_{G\cap Y}Y\, \, = \, \, \bigcap_{i=1}^t\,\bl_{G\cap V_i}V_i\, . 
\end{equation}
For the inclusion~``$\subseteq$'' note that 
$\bl_{G\cap  Y}Y\,{\subseteq}\,\bl_{G\cap V_i}V_i$ is a direct
consequence of~$Y\,{\subseteq}\, V_i$ as discussed in Step~2~(1).

For the reverse inclusion we need the following identity:
\begin{equation}\label{eq_join}
     \bigvee_{i=1}^t\, (G\wedge V_i)\,\, = \, \, G\wedge Y\, .
\end{equation}
This identity holds in any semilattice without referring to $G$ 
being an element of the building set.

Let $\alpha \,{\in}\, \cap_{i=1}^t\,\bl_{G\cap V_i}V_i$. In case
$\alpha \,{\in}\, \cap_{i=1}^t\, V_i\, {\setminus}\, (G\,{\cap}\,V_i)$,
we conclude that $\alpha\,{\in}\, Y\, {\setminus}\, (G\,{\cap}\,Y) $.
We thus assume that $\alpha$ is contained in the intersection of
exceptional divisors $\widetilde{G\,{\cap}\,V_i}$,  $i=1,\ldots,t$. 
We again switch to local considerations in
the neighborhood of a point~$y\,{\in}\,G\,{\cap}\,Y$, using that it 
carries a stratification biholomorphic to an arrangement
stratification.

Using the  description (\ref{eq_blowup}) of points in exceptional
divisors that are created by blowups in the arrangement setting, 
$\alpha\,{\in}\,\cap_{i=1}^t\,\widetilde{G\,{\cap}\,V_i}\,{\subseteq}\,
\cap_{i=1}^t\,\bl_{G\cap V_i}V_i$ means that there exist $a\,{\in}\,
\cap_{i=1}^t\, (G\,{\cap}\,V_i)$, and  
$p_i\,{\in}\, V_i\, {\setminus}\, (G\,{\cap}\,V_i)$ for $i=1,\ldots,t$,
with 
\[
\alpha \,\, = \, \, (a,\Span(p_i,G\,{\cap}\,V_i))\, \in \, 
\bl_{G\cap  V_i}V_i\, .
\] 
In particular, $\Span(p_i,G)\,{=}\,\Span(p_j,G)$ for $1\,{\leq}\,i,j\,{\leq}\, t$.
Thus,
\[
\Span(p_j,G) \, \subseteq \bigcap_{i=1}^t\, \Span(V_i,G)\, \, = \, \, \Span(Y,G)   
\] 
using the identity~(\ref{eq_join}). We conclude that there exists $y\,{\in}\,
Y\,{\setminus}\,(G\,{\cap}\,Y)$ such that 
$\Span(y,G)\, {=}\, \Span(p_j,G)$ for all $j\,{\in}\,\{1,\ldots,k\}$,
hence
\[
\alpha\, = \, (a,\Span(y,G\,{\cap}\,Y))\,\,\in\,\, \bl_{G\cap Y}Y\, .
\]

Though we are for the moment not concerned with the case of
$Y\,{\subseteq}\,G$, we note for later reference
that~(\ref{eq_intdiv})
remains true, with $\bl_Y Y\,{=}\,\emptyset$ meaning that the
intersection on the right-hand side is empty.  
Following the proof of the
inclusion~``$\supseteq$'' in~(\ref{eq_intdiv}) for
$G\,{\cap}\,Y\,{=}\,Y$,
we first find that the intersection of blowups can only contain points 
in the exceptional divisors. Assuming~$\alpha\,{\in}\,\cap_{i=1}^t\,
\widetilde{G\,{\cap}\,V_i}$ we arrive to a contradiction when
concluding that $\Span(p_j,G)\,{\subseteq}\,
\cap_{i=1}^t\, \Span(V_i,G)\,{=}\,\Span(Y,G)\,{=}\,G$ for
$j=1,\ldots,t$.

\noindent
{\bf Step 4:} {\em Any intersection of maximal strata in 
$\bl_i^{\cg}(\ca)$ occurs as an assigned space.} \newline
Every intersection involving the exceptional divisor~$\widetilde G$
occurs if we can show that all other intersections occur
(intersections that additionally involve~$\widetilde G$ then are  
assigned to corresponding elements of type~I$\!$I). 

Consider
$W\, = \, \bigcap_{i=1}^t\,\bl_{G\cap V_i}V_i$,
where the $V_i$ are maximal strata in $\bl_{i-1}^{\cg}(\ca)$; recall
here that a blowup in an empty stratum does not alter the space. 
We can assume that $\cap_{i=1}^t\, V_i\,{\not =}\, \emptyset$,
otherwise the intersection~$W$ were empty. With the
identity~(\ref{eq_intdiv}) in Step~3 we conclude that either $W\,{=}\,
\emptyset$ (in case $\cap_{i=1}^t\, V_i\,{\subseteq}\,G$) or
$W\,{=}\,\bl_{G\cap \bigcap_{i=1}^t V_i} \cap_{i=1}^t\, V_i$,
in which case it is assigned to the element $\cap_{i=1}^t\, V_i$
in $\bl_i(\cl)$.

\noindent
{\bf Step 5:} {\em $\ca_{i}^{\cg}$ is a local subspace arrangement
in $\bl_i^{\cg}(\ca)$.} 
\newline
It follows from the description~(\ref{eq_intdiv}) of strata in  
$\bl_i^{\cg}(\ca)$ that all intersections of maximal strata are
connected and smooth. It remains to show that $\ca_{i}^{\cg}$
locally looks like a subspace arrangement. Let $y\,{\in}\,\ca_{i}^{\cg}$.
We can assume that $y$ lies in the exceptional divisor~$\wti G$. 
Let $x \in G\,{\subseteq}\,\ca_{i-1}^{\cg}$ be the image of~$y$ under 
the blowdown map.

We first give a local description around $x$ in $\ca_{i-1}^{\cg}$.
By induction hypothesis, there exists a neighborhood~$N$ of~$x$, and 
an arrangement of linear subspaces~$\cb$ in~$\dc^n$
such that the pair $(N,\ca_{i-1}^{\cg}\,{\cap}\,N)$ is biholomorphic
to the pair~$(\dc^n,\cb)$. We can assume that under this biholomorphic 
map, $x$ is mapped to the origin. Let $T\,{=}\,\bigcap_{B\in \cb}\,B$ 
and note that $G\,{\cap}\,N$ is mapped to some subspace~$\ga$ in $\cb$.

With $G$ being maximal in the building set for~$\ca_{i-1}^{\cg}$,
$\cb/T$ is a product arrangement with one of the factors being an
arrangement in~$\ga/T$. More precisely, there exists a subspace 
$\ga'\,{\subseteq}\,\dc^n$, and two subspace arrangements, $\cc$ in $\ga/T$
and $\cc'$ in $\ga'/T$, such that
\begin{itemize}
\item[(1)] $\ga/T\, \oplus\, \ga'/T\, \oplus\, T\, = \, \dc^n$,  
\item[(2)] $\cb\, =\ 
\{A\, \oplus\, \ga'/T\, \oplus\, T\,| \, A\,{\in}\,\cc\} \, \cup \,
\{\ga/T\, \oplus\, A'\, \oplus\, T\,| \, A'\,{\in}\,\cc'\}$. 
\end{itemize}

Blowing up $G$ in $\bl_{i-1}^{\cg}(\ca)$ locally corresponds to
blowing up $\Gamma$ in $\dc^n$. Let $t$ be the point on the special divisor
$\wti \ga$ corresponding to~$y \in \wti G$, thus $t$ maps to the origin 
in~$\dc^n$ under the blowdown map. A neighborhood of~$t$ in $\bl_{\ga}\dc^n$ 
is an $n$-dimensional open ball which can be parameterized as a direct sum
\[
M \,\oplus \, M' \, \oplus \, I \, \oplus \,T \, .
\]
Here, $M$ is an open ball around~$0$ in $\ga/T$, $M'$ is an open ball on the
unit sphere in  $\ga'/T$ around the point of intersection with the 
line~$\langle p \rangle$
in  $\ga'/T$ that defines $t$ as a point in the exceptional divisor, 
$t=(0, \Span(p,\ga)) \in \wti \ga$ (compare (\ref{eq_blowup})), 
and $I$ an open unit ball in~$\dc$.

The maximal strata in this neighborhood are the following:
\begin{itemize}
\item[$\circ$]
      the hyperplane $M \,\oplus \, M' \, \oplus \, \{0\} \, \oplus \,T$,
      as the exceptional divisor,
\item[$\circ$]
      $(M \cap A)\,\oplus \, M' \, \oplus \, I \, \oplus \,T$, replacing
      $A\, \oplus\, \ga'/T\, \oplus\, T$ after blowup,
\item[$\circ$]
      $M \,\oplus \, (M'\cap A') \, \oplus \, I \, \oplus \,T$, replacing
      $\ga/T\, \oplus\, A'\, \oplus\, T$ after blowup for $A'\neq 0$.
\end{itemize} 

This proves that around $t$ in $\bl_{\ga}\dc^n$ we have the structure of a 
local subspace arrangement, which in turn shows the local arrangement 
property around $y$ in $\ca_i^{\cg}$.

\noindent
{\bf Step 6:} {\em $\cg$ is a building set for $\ca_i^{\cg}$
 in the sense of \/{\rm Definition~\ref{lsabsdf}}.}\newline
$\cg$ is a~combinatorial building set by 
Proposition~\ref{prop_singleblow}. Complementing this with the 
dimension information about the strata, we conclude, by
Proposition~\ref{combgeomprop}(2), that $\cg$ is a~geometric 
building set.
\qed


\subsection{Simplicial Resolutions Of Toric Varieties}
\label{ssect_tv}
$\,$
\vskip3pt

The study of toric varieties has proved to be a field of fruitful
interplay between algebraic and convex geometry: toric
varieties are determined by rational polyhedral fans, and many of their 
algebraic geometric properties are reflected by combinatorial
properties of their defining fans. 

We recall one such correspondence -- between subdivisions of fans and
special toric morphisms -- and show that so-called stellar subdivisions 
are instances of combinatorial blowups. This allows us to apply our
Main Theorem in the present context: Given a polyhedral fan, we
specify a class of {\em simplicial\/} subdivisions, and, interpreting
our notions of building sets and nested sets, we describe the 
incidence combinatorics of the subdivisions in terms of the
combinatorics of the initial fan. For background material on
toric varieties we refer to the standard sources~\cite{Da,Od,Fu,Ew}.

Let $X_{\Sigma}$ be a toric variety defined by a rational polyhedral
fan~$\Sigma$. Any subdivision of~$\Sigma$ gives rise to a proper,
birational toric morphism between the associated toric 
varieties~(cf~\cite[5.5.1]{Da}). In particular, simplicial
subdivisions yield toric morphisms from quasi-smooth toric varieties
to the initial variety -- so-called {\em simplicial resolutions}.
Quasi-smooth toric varieties being rational homology manifolds, such
morphisms can replace smooth resolutions for (co)homological
considerations.

We define a particular, elementary, type of subdivisions:

\begin{df} \label{def_stellsd}
Let $\Sigma\,{=}\,\{\sigma\}_{\sigma \in \Sigma}\,{\subseteq}\,\dr^d$
be a polyhedral fan, i.e., a collection of closed polyhedral
cones~$\sigma$ in~$\dr^d$ such that $\sigma\,{\cap}\,\tau$ is a cone
in~$\Sigma$ for any $\sigma, \tau\,{\in}\,\Sigma$. Let ${\rm cone}(x)$ be
a ray in~$\dr^d$ generated by $x\,{\in}\,{\rm relint}\, \tau$ for some 
$\tau\,{\in}\,\Sigma$. The {\em stellar subdivision\/} ${\rm sd}(\Sigma,x)$
of $\Sigma$ in $x$ is given by the collection of cones
\[
(\,\Sigma\,\setminus\, {\rm star}(\tau, \Sigma)\, )
\,\, \cup \, \, 
\{\, {\rm cone}(x,\rho)\, | \, \rho \subseteq \sigma \,\,
                               \mbox{ for some } \,\,
                               \sigma \in {\rm star}(\tau, \Sigma)\,    
\}\, ,
\] 
where ${\rm star}(\tau, \Sigma)\,{=}\,\{\sigma\,{\in}\,\Sigma\, | \, 
\tau\,{\subseteq}\,\sigma\}$, and ${\rm cone}(x,\rho)$ the closed polyhedral 
cone
spanned by~$\rho$ and~$x$. If only concerned with the combinatorics of 
the subdivided fan, we also talk about stellar subdivision of $\Sigma$ 
in $\tau$, ${\rm sd}(\Sigma,\tau)$, meaning any stellar subdivision in~$x$
for~$x\,{\in}\,{\rm relint}\, \tau$. 
\end{df} 

\begin{prop}
Let $\cf(\Sigma)$ be the face poset of a polyhedral fan~$\Sigma$,
i.e., the set of closed cones  in~$\Sigma$ ordered by inclusion, together with the
zero cone~$\{0\}$ as a minimal element. For $\tau\,{\in}\,\Sigma$,
the face poset of the stellar subdivision of $\Sigma$ in $\tau$ can be
described as the combinatorial blowup of $\cf(\Sigma)$ at~$\tau$:
\[
  \cf({\rm sd}(\Sigma, \tau))\,\, = \,\, \bl_{\tau}\cf(\Sigma)\, . 
\] 
\end{prop}

\pr
Removing ${\rm star}(\tau, \Sigma)$ from~$\Sigma$ corresponds to
removing $\cf(\Sigma)_{\geq \tau}$ from $\cf(\Sigma)$, adding cones as 
described in Definition~\ref{def_stellsd} corresponds to extending 
$\cf(\Sigma)\,{\setminus}\,\cf(\Sigma)_{\geq \tau}$ by
elements~$[\tau,\rho]$ for $\rho\,{\in}\,\cf(\Sigma)$, 
$\rho\,{\subseteq}\,\sigma$  for some $\sigma \,{\in}\,{\rm star}(\tau, \Sigma)$.
The comparison of order relations is straightforward. \qed

\vspace{0.2cm}
We apply our Main Theorem to the present context.

\begin{thm}
Let~$\Sigma$ be a polyhedral fan in~$\dr^d$ with 
face poset~$\cf(\Sigma)$. Let $\cg\,{\subseteq}\,\cf(\Sigma)$ be a building 
set 
of~$\cf(\Sigma)$ in the sense of Definition~\ref{df_buildg}, $\cn(\cg)$ the
complex of nested sets in~$\cg$ (cf.\ Definition~\ref{df_nested}).
Then, the consecutive application of stellar subdivisions in every cone 
$G\,{\in}\,\cg$ in a non-increasing order yields a simplicial subdivision
of~$\Sigma$ with face poset equal to the face poset of~$\cn(\cg)$.
\end{thm}

\noindent
As examples of building sets for face lattices of 
polyhedral fans let us mention:
\begin{itemize}
\item[(1)] the full set of faces, with the corresponding complex of nested sets 
           being the order complex of~$\cf(\Sigma)$ (stellar subdivision in all 
           cones results in the barycentric subdivision of the fan);
\item[(2)] the set of rays together with the non-simplicial faces of~$\Sigma$;
\item[(3)] the set of irreducible elements in $\cf(\Sigma)$: the set of rays 
           together with all faces of~$\Sigma$ that are not products 
           of some of their proper faces.
\end{itemize}

\begin{remrm} \label {rm_smooth_tv}
For a smooth toric variety~$X_{\Sigma}$, the union of closed codimension~1 
torus orbits is a local subspace arrangement, in particular, the codimension~1 
orbits form a divisor with normal crossings, \cite[p.~100]{Fu}. 
The intersection stratification of this local arrangement coincides
with the torus orbit stratification of the toric variety. For any face 
$\tau$ in the defining fan~$\Sigma$, the torus orbit $\co_{\tau}$
together with all orbits corresponding to rays in~$\Sigma$ form 
a~geometric building set. Our proof in~\ref{ssect_tracg_inc} applies
in this context with $\co_{\tau}$ playing the role of~$G$. We conclude 
that under blowup of~$X_{\Sigma}$ in the closed torus orbit
$\co_{\tau}$, the incidence combinatorics of torus orbits changes
exactly in the way described by a stellar subdivision of~$\Sigma$ 
in~$\tau$. This is the combinatorial part of the well-known fact that
in the smooth case, the blowup of $X_{\Sigma}$ in a torus orbit 
$\co_{\tau}$ corresponds to a regular stellar subdivision of the 
fan~$\Sigma$ in $\tau$~\cite{MO}.
 \end{remrm}


\section{An outlook}
\label{sect5}

\subsection{Models for real subspace arrangements and stratified manifolds}
\label{ssectGaiffi}
\mbox{ }

In the spirit of the De Concini-Procesi wonderful model construction
for subspace arrangements, Gaiffi~\cite{Ga2} presents a model
construction for the complement of arrangements of real linear subspaces
 modulo $\mathbb R^+$: Given a central subspace
arrangement $\ca$ in some Euclidean vector space~$V$, denote by
$\widehat \cm(\ca)$ the quotient of its complement by $\mathbb R^+$.
Denote the unit sphere in $V$ by $S(V)$, and consider, for a given
(geometric) building set $\cg$ in $\cl(\ca)$, the embedding
\[
   \rho\,: \, \, \widehat \cm(\ca) \, \, \longrightarrow \, \, 
                            S(V) \,\times \prod_{G\in \cg}\, G\cap S(V)\, .
\]
The map is obtained by composing the natural section $\widehat
\cm(\ca)\, \rightarrow \, \cm(\ca)$, $[x] \mapsto \frac{x}{|x|}$, with
a projection onto each factor of the right-hand side product. Denote
the closure of this map by $Y_{\cg}$. $Y_{\cg}$ is shown to be a
manifold with corners, which enjoys much of the properties familiar
from the projective setting: the boundary of $Y_{\cg}$ is stratified
by codimension~$1$ manifolds with corners indexed with building set
elements and having non-empty intersection whenever the index set is
nested with respect to~$\cg$. The set-up allows for a straightforward
generalization to mixed subspace and halfspace arrangements motivated
by compactifications of configuration spaces in work of
Kontsevich~\cite{Ko}. A step aside from classical (linear)
arrangements, our combinatorial framework still applies is this
context.

In a second part of his paper, Gaiffi extends the previous
construction to conically stratified manifolds with corners. Replacing
the explicit construction of taking the closure of an embedding into a
product of spheres as above, he describes a sequence of ``real
blowups'' in the sense of Kuperberg~\& Thurston~\cite{KuTh}. The
sequence is prescribed by the choice of a subset of strata in the
original manifold that is a combinatorial building set in our sense.
The resulting space is a manifold with corners with its boundary
stratified by codimension~$1$ manifolds with corners that are indexed
by the building set elements, and intersections being non-empty if and
only if the corresponding index sets are nested.

\subsection{A graded algebra associated with a finite lattice}
\mbox{ }

In a joint paper of Yuzvinsky and the first author~\cite{FY}, we start out
from the combinatorial notions of building sets and nested sets given in
the present paper and define a commutative graded algebra in purely 
combinatorial terms:

\begin{df}
  For a finite lattice $\cl$, $\mathfrak A$ its set of atoms, and
  $\cg$ a combinatorial building set in $\cl$, define the algebra
  $D(\cl,\cg)$ as the quotient of a polynomial algebra over $\mathbb
  Z$ with generators in $1$-$1$ correspondence with the elements
  of~$\cg$:
\[ 
  D(\cl,\cg) \, \, := \, \, 
             \mathbb Z\,[\{x_{G}\}_{G\in \cg}] \, \Big/ \, \ci \, ,
\]
where the ideal of relations $\ci$ is generated by 
\begin{eqnarray*}
    \prod_{i=1}^t \,x_{G_i}\, , & & 
             \mbox{for }\,\,\{G_1,\ldots ,G_t\} \mbox{ not nested}\,,\\
    \sum_{G\geq H}\, x_G\, , & &  
                            \mbox{for }\, H \in \mathfrak A \, .  
\end{eqnarray*}
\end{df}

For $\cl$ the intersection lattice of an arrangement of complex
hyperplanes~$\ca$ and~$\cg$ its minimal building set, this algebra was shown
to be isomorphic to the integer cohomology algebra of the compact wonderful
arrangement model in~\cite[1.1]{DP2}. We show in~\cite{FY} that the  
algebra in fact is isomorphic to the cohomology algebra of the arrangement 
model for {\em any\/} choice of a building set in the intersection lattice.

Going beyond the arrangement context, we can provide yet another geometric 
interpretation of the algebras~$D(\cl,\cg)$: For an arbitrary atomic 
lattice and a given combinatorial building set we construct a smooth, 
non-compact toric variety $X_{\Sigma(\cl,\cg)}$ and show that its Chow 
ring is isomorphic to the algebra~$D(\cl,\cg)$. 

In a sense, this is a prototype result of what we had hoped for when working 
on our combinatorial framework: to provide the outset
for going beyond the geometric context of resolutions and yet get back
to it in a different, elucidating, and, other than via the abstract 
combinatorial detour, seemingly unrelated way.


\end{document}